\documentclass[11pt,a4paper]{article}

\usepackage[english]{babel}
\usepackage{amsmath,amssymb,amsfonts,a4}
\usepackage{eucal,bbm,stmaryrd}
\usepackage{hyperref}
\usepackage{xcolor}
\usepackage{tikz}

\begin{document}

\numberwithin{equation}{section}

\newtheorem{Thm}{Theorem}
\newtheorem{Prop}{Proposition}
\newtheorem{Def}{Definition}
\newtheorem{Lem}{Lemma}
\newtheorem{Rem}{Remark}
\newtheorem{Cor}{Corollary}
\newtheorem{Con}{Conjecture}

\newcommand{\Thmautorefname}{Theorem}
\newcommand{\Propautorefname}{Proposition}
\newcommand{\Defautorefname}{Definition}
\newcommand{\Lemautorefname}{Lemma}
\newcommand{\Remautorefname}{Remark}
\newcommand{\Corautorefname}{Corollary}
\newcommand{\Conautorefname}{Conjecture}

\newcommand{\Pf}{\noindent{\bf Proof: }}
\newcommand{\qed}{\hspace*{3em} \hfill{$\square$}}

\newcommand{\N}{\mathbbm{N}}
\newcommand{\Z}{\mathbbm{Z}}
\newcommand{\R}{\mathbbm{R}}

\newcommand{\E}{\mathbbm{E}}
\renewcommand{\P}{\mathbbm{P}}

\newcommand{\Bin}{\mathit{Bin}}

\newcommand{\G}{\Gamma}

\newcommand{\la}{\lambda}
\newcommand{\La}{\Lambda}
\newcommand{\si}{\sigma}
\newcommand{\al}{\alpha}
\newcommand{\be}{\beta}
\newcommand{\ep}{\epsilon}
\newcommand{\ga}{\gamma}
\newcommand{\de}{\delta}
\newcommand{\De}{\Delta}
\newcommand{\ph}{\varphi}
\newcommand{\om}{\omega}
\renewcommand{\th}{\theta}
\newcommand{\vth}{\vartheta}
\newcommand{\ka}{\kappa}

\newcommand{\F}{\mathcal{F}}
\newcommand{\cZ}{\mathcal{Z}}
\newcommand{\cT}{\mathcal{T}}
\newcommand{\cX}{\mathcal{X}}
\newcommand{\cY}{\mathcal{Y}}
\newcommand{\cW}{\mathcal{W}}
\newcommand{\tcW}{\tilde{\cW}}
\newcommand{\cV}{\mathcal{V}}
\newcommand{\cC}{\mathcal{C}}
\newcommand{\cN}{\mathcal{N}}
\newcommand{\da}{\dagger}
\newcommand{\Mda}{M^{\dagger}}
\newcommand{\tM}{\tilde{M}}
\newcommand{\xda}{x_\da}
\newcommand{\xst}{x_*}

\newcommand{\Aut}{Aut}

\newcommand{\Nm}{N_{\wedge}}

\newcommand{\stm}{\setminus}
\newcommand{\lra}{\leftrightarrow}
\newcommand{\lraa}[1]{\stackrel{#1}{\lra}}
\newcommand{\Ra}{\Rightarrow}
\newcommand{\Lra}{\Leftrightarrow}

\newcommand{\cp}{\, \square \,}

\thispagestyle{plain}
\title{Comparing the number of infected vertices\\ 
in two symmetric sets for Bernoulli percolation\\
(and other random partitions)}
\author{Thomas Richthammer \footnote{Institut f\"ur Mathematik, Universit\"at Paderborn}}
\maketitle
\begin{abstract}
For Bernoulli percolation on a given graph $G = (V,E)$
we consider the cluster of some fixed vertex $o \in V$.
We aim at comparing the number of vertices of this cluster in the set $V_+$ and in the set $V_-$, where $V_+,V_- \subset V$
have the same size. 
Intuitively, if $V_-$ is further away from $o$ than $V_+$, 
it should contain fewer vertices of the cluster. 
We prove such a result in terms of stochastic domination, 
provided that $o \in V_+$, and $V_+,V_-$ satisfy some strong symmetry conditions, 
and we give applications of this result
in case $G$ is a bunkbed graph, a layered graph, the 2D square lattice or a hypercube graph. 
Our result extends to fairly general random partitions, e.g. as induced by Bernoulli site percolation or the random cluster model.

\bigskip 

\end{abstract}



\section{Introduction}

We consider Bernoulli (bond) percolation on a graph $G = (V,E)$ with parameter $p \in (0,1)$. This stochastic process is given 
in terms of a family $(\cZ_e)_{e \in E}$ of independent $\{0,1\}$-valued random variables such that $\P_p(\cZ_e = 1) = p$. 
Here the index in the notation for probabilities indicates the $p$-dependence. We may think of percolation as the random subgraph of $G$ with vertex set $V$ and edge set $E_\cZ = \{e \in E: \cZ_e = 1\}$. The connectedness relation in this random subgraph usually is denoted by $\lra$, and the connected component of a vertex $v \in V$ is called the cluster of $v$ and denoted by $\cC_v$. Thus 
$$
\cC_v = \{u \in V \!: \exists n \ge 0, v_0,...,v_n \in V \!: v_0 = v, v_n = u, 
\forall 0 \le i < n \!: \cZ_{v_iv_{i+1}} = 1\}.  
$$
In the case that $G$ is an infinite graph with some regularity, e.g. $G = \Z^d$, Bernoulli percolation is very well studied  and there still are various interesting open problems (e.g. see \cite{G1}). 
We will consider the cluster $\cC_o$ of a fixed vertex $o \in V$. 
$o$ could be thought of as the origin of an infection that can be transmitted via edges $e \in E$ such that $\cZ_e = 1$, so that 
$\cC_o$ represents the set of infected vertices. 
In case of $G$ being a regular lattice such as $\Z^d$ 
there are classical results on the size of $\cC_o$ (e.g. see \cite{G1}) 
and the asymptotic shape of $\cC_o$ (e.g. see \cite{ACC}). 
Our focus is on monotonicity properties of the shape. 

To motivate our interest, let us consider the cluster 
of $o = (0,0)$ for Bernoulli percolation on $G = \Z^2$, 
and let $\cW_n$ denote  the number of infected vertices 
in the layer $V_n := \{n \} \times \Z$, i.e. $\cW_n := |\cC_o \cap V_n|$.  
Following the intuition that vertices closer to $o$ can be infected more easily than vertices further away, 
it is reasonable to expect that $\cW_n$ decreases with the distance $|n|$ of the layer to the origin $o$. 
Of course, one has to make sense of the comparison of the number of infected vertices in different layers: 
One might compare average values and thus expect that 
$\E_p(\cW_m) \le \E_p(\cW_n)$ for all $0 \le n \le m$. 
Alternatively one might compare the preference of the distributions to take greater values (in the sense of stochastic domination) and thus even expect that $\cW_m \preceq \cW_n$ w.r.t. $\P_p$ for all $0 \le n \le m$, which is defined to mean $\E_p(f(\cW_m)) \le \E_p(f(\cW_n))$ for all increasing, bounded functions $f: \{0,1,2,...\} \cup \{\infty\} \to \R$. 
It turns out that both of these conjectures are open problems.  
There is an obvious link to the monotonicity of connection probabilities. E.g. in the above example it should be expected 
that $\P_p(o \lra (n,0))$ is decreasing in $n$ for $n \ge 0$.
So far this is only known if $p$ is sufficiently small,   
see  \cite{LPS} and the references therein. 
Another famous example of this kind of monotonicity property 
is the bunkbed conjecture, 
which will be discussed in the following section. 
It seems that these kind of monotonicity properties in general 
are difficult to establish (if they hold at all), 
which at first glance might seem somewhat surprising. 
The contribution of this paper is the comparison of the infected number of vertices in two subsets $V_+,V_- \subset V$ in the special case that $o$ is contained in $V_+$ and assuming that $V_+,V_-$ are highly symmetric. 
Among other interesting applications of our result we obtain that 
$\E_p(\cW_m) \le \E_p(\cW_0)$ for all $m \ge 0$ in the above notation 
for Bernoulli percolation on $G = \Z^2$. 

We will state our general result in the following section and provide 
several applications, namely in case that $G$ is a bunkbed graph, a layered graph, the 2D square lattice or a hypercube graph. 
While our interest is mainly in percolation, the proof of our main result does not use any specific properties of this model, but relies on general probabilistic techniques such as conditioning, a symmetrization procedure, and a combinatorial argument concerning group actions on pairs of sets. 
Thus in Section \ref{Sec:action} we will briefly revisit some aspects of actions of groups on sets and state and proof our combinatorial result. In Section \ref{Sec:partitions} we will 
state and prove our generalized result on random partitions that allows for an application not only in Bernoulli bond percolation, 
but in many other processes that produce random partitions such as Bernoulli site percolation or the random cluster model (see \cite{G2}). In Section \ref{Sec:percolation} we will return 
to the setting of Bernoulli percolation and infer the results 
presented in Section~\ref{Sec:results} from our generalized results.

\newpage 

\section{Results} \label{Sec:results}

Let us start by reviewing concepts and notations concerning graphs 
that will be used throughout the paper.  
Let $G = (V,E)$ be a fixed (simple) graph. 
We will assume $G$ to be connected and locally finite, 
i.e. for every $v \in V$ the degree $\deg(v)$ is finite. 
We use the usual notation for edges, e.g. $vw := \{v,w\}$ 
for the edge connecting $v,w \in V$. 
We let $d(v,w)$ denote the graph distance of $v,w \in V$, i.e. 
the number of edges  of the shortest path connecting $v$ and $w$. 
The inner boundary of a set $A \subset V$ of thickness $k \in \N$ will be denoted by
$$
\partial_k A := \{v \in A: \exists w \in A^c: d(v,w) \le k\}
$$ 
The set of graph automorphisms $\Aut(G)$ is the set 
of all bijections on $V$ such that the edges are preserved in that 
for all $v,w \in V$ we have $vw \in E$ iff $\ph(v)\ph(w) \in E$. 
A set $\G \subset \Aut(G)$ 
is said to act transitively on a set $V' \subset V$ iff
$$
\forall \ph \in \G: \ph(V') = V' \quad \text{ and } \quad 
\forall v,w \in V' \exists \ph \in \G: \ph(v) = w. 
$$
Let 
\vspace*{ - 0.3 cm}
$$
\G_v(w) = \{\ph(w): \ph \in \G \text{ s.t. } \ph(v) = v\}
$$
denote the orbit of vertex $w \in V$ under all graph automorphisms in $\G$ fixing $v \in V$. 
We are now ready to state the main result of our paper: 

\begin{Thm} \label{Thm:Graphs} 
We consider Bernoulli percolation with parameter $p \in (0,1)$ on a connected, locally finite graph $G = (V,E)$.  
Let $\cC := \cC_o$ denote the cluster of a fixed vertex 
$o \in V$. 
Let $V_+,V_- \subset V$ be two disjoint sets of vertices 
and let $\G$ be a subgroup of $\Aut(G)$ such that we have 
\begin{itemize}
\item[$(\Gamma 1)$] $\forall \ph \in \G: \ph(V_+),\ph(V_-) \in \{V_+,V_-\}$,
\item[$(\Gamma 2)$] $\G$ acts transitively on $V_\pm := V_+ \cup V_-$ and 
\item[$(\Gamma 3)$] $\forall v \in V_+,w \in V_-: |\G_v(w)| = |\G_w(v)|$
\end{itemize}
Let $\cC_+ := \cC \cap V_+$, $\cC_- := \cC \cap V_-$ and 
$\cC_\pm := \cC_+ \cup \cC_-$. 
If $o \in V_+$, then  
we have 
\begin{align*}
&|\cC_-| \preceq |\cC_+| \quad \text{  w.r.t. } 
\P_p(\,.\,|\,|\cC_\pm| < \infty) \quad \text{ and in particular}\\
&\E_p(\cC_-|\,|\cC_\pm| < \infty) \le \E_p(\cC_+|\,|\cC_\pm| < \infty)\quad \text{ and } \quad \E_p(\cC_-) \le \E_p(\cC_+).
\end{align*}
\end{Thm}

\begin{Rem} Some comments on the above theorem: 
\begin{itemize}
\item 
The sets $V_+,V_-$ in the theorem are symmetric 
to each other in that there is a graph automorphism $\ph \in \G$ such that $\ph(V_+) = V_-$ (by $(\G 1)$ and $(\G 2)$), 
and $V_+$ ist closer to $o$ than $V_-$ by the (rather strong) assumption that $o \in V_+$. 
Thus the cluster $\cC$ of $o$ could be expected to contain more vertices of $V_+$ than vertices of $V_-$. 
The theorem gives a precise statement of this comparison in terms of stochastic domination and a comparison of expectations. 
It also spells out a precise symmetry assumption on $V_+, V_-$ 
in terms of the existence of a subgroup $\G \subset \Aut(G)$ 
with certain properties. 
\item 
In the conclusion of the theorem both  
conditional and unconditional expectations are compared. 
In case of infinite sets $V_+,V_-$, 
in the supercritical phase of percolation only the former is interesting, 
while in the subcritical phase there is no difference between them. 
\item 
Given $V_+,V_-$ it usually makes sense to choose $\G$ as large as possible, 
i.e. $\G := \{\ph \in \Aut(G): \ph(V_+),\ph(V_-) \in \{V_+,V_-\}\}$. 
Then $(\G 1)$ holds trivially and we have the best chance for $(\G 2)$ to hold. 
Indeed, the only advantage for a different choice for $\G$ 
is that it may be easier to verify $(\G 3)$. 
\item 
$(\G 1)$ and $(\G 2)$ are strong assumptions on the internal 
symmetries of $V_+$ and $V_-$:   
They imply that $\G' = \{\ph \in \G: \ph(V_+) = V_+, \ph(V_-) = V_-\}$ acts transitively on $V_+$ and on  $V_-$. 
Nevertheless we will see that this still allows for interesting applications. 
\item 
$(\G 3)$ is a weak assumption. 
While it is possible to engineer graphs and sets $V_+,V_-$ so that for all choices of $\G$
satisfying $(\G 1)$ and $(\G 2)$ we do not have $(\G 3)$, 
in most applications $(\G 3)$ is easily verified. 
In the following proposition we collect some sufficient conditions for $(\G 3)$.  
%
%
%
\item 
As already remarked in the introduction, 
our result extends to rather general random partitions 
of $V_\pm = V_+ \cup V_-$, see
Theorem \ref{Thm:partitions} in Section \ref{Sec:partitions}. 
\end{itemize}
\end{Rem}

\begin{Prop} \label{Prop:orbstab3}
Let $G = (V,E)$ be a connected, locally finite graph, let 
$V_\pm \subset V$ and let $\G$ be a subgroup of $\Aut(G)$ acting transitively on $V_\pm$. 
If one of the following conditions is satisfied, 
then $\forall v,w \in V_\pm : |\G_v(w)| = |\G_w(v)|$. 
\begin{itemize}
\item[$(\G 3a)$] $\forall v,w \in V_\pm \exists \ph \in \G\!: \ph(v) = w, \ph(w) = v.$
\vspace*{ - 0.2 cm}
\item[$(\G 3b)$] $\displaystyle \forall k \ge 1 : 
 \inf\Big\{ \frac{|V_\pm \cap \partial_k K|}{|V_\pm \cap K|} : K \subset V \text{ s.t. } |K| < \infty, V_\pm \cap K \neq \emptyset  \Big\}  =   0$.
\end{itemize}
\end{Prop}
In the following we give applications of Theorem~\ref{Thm:Graphs}. 
We concentrate on graph classes that appear elsewhere in the literature with some connection to monotonicity questions. 
These graphs turn out to be Cartesian products:   

\begin{Def}
For graphs $G_1 = (V_1,E_1),...,G_n = (V_n,E_n)$ 
their Cartesian product $G_1 \cp  ... \cp G_n$ is the graph $(V,E)$ with $V = V_1 \times ... \times  V_n$ and 
$E = \{(v_1,...,v_n)(w_1,...,w_n)\!: \exists i \in \{1,...,n\}\!: 
v_iw_i \in E_i, \forall j \neq i \!: v_j = w_j \in V_j\}$.
\end{Def}
We note that any choice of $\ph_i \in \Aut(G_i)$ $(1 \le i \le n)$
gives an automorphism $(\ph_1,...,\ph_n)$ on the product graph. 
In the first two applications we use the following graph properties 
(related to ($\G$3a) and  ($\G$3b) from Proposition \ref{Prop:orbstab3}):  

\begin{Def}
A graph $G = (V,E)$ is called 
\begin{itemize}
\item 
vertex-transitive 
iff $\forall v,w \in V \exists \ph \in \Aut(G): \ph(v) = w$. 
\item  
vertex-swap-transitive
iff $\forall v,w \in V \exists \ph \in \Aut(G): \ph(v) = w, \ph(w) = v$. 
\vspace*{-0.2 cm}
\item 
$*$-amenable, iff $\displaystyle \forall k \ge 1 : \inf \Big\{ \frac{|\partial_k K|}{|K|} : K \subset V \text{ s.t. }
0 < |K| < \infty \Big\} = 0$. 
\end{itemize}
\end{Def}
Vertex-transitivity is a well known property, 
vertex-swap-transitivity is a more restrictive variant of that, and $*$-amenability is a variant of the usual amenability property. 
We note that every finite graph is $*$-amenable
since $\partial_k V = \emptyset$. 

\begin{Cor} \label{Cor:bunkbed} Bunkbed graphs. 
Let $G = (V,E)$ be a connected, locally finite graph, 
and $L_2$ the line graph on two vertices, 
i.e. $L_2 = (\{0,1\}, \{\{0,1\}\})$.
We consider the bunkbed graph $G \cp L_2$ and  
$o = (o',0)$ for some $o' \in V$. 
If $G$ is either a $*$-amenable vertex-transitive graph
or a vertex-swap-transitive graph, 
then the conclusion of Theorem \ref{Thm:Graphs} holds for 
$V_+ = V \times \{0\}$ and $V_- = V \times \{1\}$. 
\end{Cor}
In bunkbed graphs one might expect a stronger monotonicity property: 
For a connected and locally finite graph 
$G = (V,E)$, 
all $p \in (0,1)$ and arbitrary $o' \in V$ 
the corresponding Bernoulli percolation on 
$G \cp L_2$ might satisfy
$$
\forall v \in V: \P_p((o',0) \lra (v,0)) \ge \P_p((o',0) \lra (v,1)). 
$$
This property immediately implies $\E_p(\cC_+) \ge \E_p(\cC_-)$ for $V_+ = V \times \{0\}$ and $V_- = V \times \{1\}$, since $\cC_+ = \sum_{v \in V}
1_{\{o \lra (v,0)\}}$ and similarly for $\cC_-$.
It was conjectured that this stronger monotonicity property holds for any (finite) graph. 
While it can be shown to hold for a few simple classes of graphs 
(see \cite{HL} or \cite{R} and the references therein),  
recently the conjecture has been disproven by carefully constructing a finite graph $G$ 
such that the above monotonicity does not hold (see  \cite{GPZ}, 
building on a related counterexample in \cite{H}). 
We note that at this point it is unclear, 
whether $\E_p(\cC_+) \ge \E_p(\cC_-)$ for all (finite) graphs $G$, 
but the counterexample to the bunkbed conjecture gives some reason 
to doubt that. 
In this context even the comparison of expectations 
from Corollary \ref{Cor:bunkbed} might be considered an 
interesting result. For the comparison in terms of stochastic domination there is no immediate relation to the 
comparison of connection probabilities. 

%
%
\begin{Cor} \label{Cor:layer} Layered graphs. 
Let $G = (V,E)$ be a connected, locally finite graph, 
and $L_\infty$ the bi-infinite line graph, 
i.e. $L_\infty = (\Z, \{\{i,i+1\}: i \in \Z\})$.
We consider the layered graph $G \cp L_\infty$ and  
$o = (o',0)$ for some $o' \in V$. 
If $G$ is either a $*$-amenable, vertex transitive graph or 
a vertex-swap-transitive graph, then the conclusion of Theorem \ref{Thm:Graphs} holds for the following choices of $V_+,V_-$: 
\begin{itemize}
\item[(a)] 
$V_+ = V \times \{0\}$ and $V_- = V \times \{m\}$ for $m \neq 0$, 
\item[(b)] $V_+ = V \times n\Z$ and $V_- = V \times (n\Z + m)$ for $1 \le m < n$. 
\item[(c)] 
$V_+ = V \times (2n\Z \cup (2n\Z + m))$ and 
$V_- = V \times ((2n\Z + n) \cup (2n\Z +n+ m))$ for 
$0 < m < 2n$ with $m \neq n$.
\end{itemize}
\end{Cor}
As with the bunkbed conjecture, for many interesting finite graphs $G = (V,E)$ (and maybe even for all finite graphs), for all $p \in (0,1)$ and arbitrary $o' \in V$ the corresponding Bernoulli percolation on the layered graph $G \cp L_\infty$ might satisfy 
\begin{align*}
&\forall v \in V \forall m > n \ge 0: \P_p(o \lra (v,n)) \ge \P_p(o \lra (v,m)) \quad \text{ and thus also }\\
&\forall m > n \ge 0: \E_p(|\cC_o \cap  V \times \{n\}|) \ge \E_p(|\cC_o \cap  V \times \{m\}|).
\end{align*}
The latter would immediately imply the expectation part of the above corollary in case of (a). 
For some related monotonicity results for layered graphs 
we refer to \cite{KR}. 
Also we would like to point out that the special case of $G$ 
being a regular tree (as e.g. considered in \cite{GN}) 
is covered in the above Corollary, since regular trees are vertex-swap-transitive. 

Corollaries \ref{Cor:bunkbed} and \ref{Cor:layer}
can easily be generalized, e.g. considering graphs with  symmetric layers, 
where the layers are connected by edges in a different (but still sufficiently symmetric) way, 
or considering general product graphs $G = G_1 \cp G_2$ 
or indeed $G = G_1 \cp ... \cp G_n$, where 
the $G_i$ satisfy suitable conditions. 
Instead of pursuing possible gerneralizations, 
we now turn to a concrete graph, namely the graph that is most prominently featured in Bernoulli percolation.

\begin{Cor} \label{Cor:Z2} 2D square lattice $\Z^2$. 
We consider $\Z^2 = L_\infty \cp L_\infty$ and $o = (0,0)$. The conclusion of Theorem \ref{Thm:Graphs} holds for 
\begin{itemize}
\item[(a)] $V_+ = \Z v$, $V_- = w + \Z v$ for arbitrary 
$v \neq o, w \notin \Z v$. 
\item[(b)] $V_+ \!= \Z v + \Z u, V_- \!=  w + V_+$, for 
lin. independent $v,u \in \Z^2$ and $w \notin V_+$.  
\end{itemize}
\end{Cor}
\begin{figure}[htb!] 

\centering
\vspace*{0.01 cm}
\begin{tikzpicture}[scale = 0.128]
\clip(-2.5,-2.5) rectangle (21.5,11.5);
\draw[lightgray, very thin] (-6,-6) grid (22,12) ;
\draw[fill] (0,0) circle (12 pt); 
\draw[fill] (6,3) circle (12 pt); 
\draw (-1,2) circle (12 pt); 
\draw (5,5) circle (12 pt); 
\end{tikzpicture}
\begin{tikzpicture}[scale = 0.128]
\clip(-2.5,-2.5) rectangle (21.5,11.5);
\draw[lightgray, very thin] (-6,-6) grid (22,12) ;
\draw[fill] (0,0) circle (12 pt); 
\draw[fill] (5,5) circle (12 pt); 
\draw (-1,2) circle (12 pt); 
\draw (6,3) circle (12 pt); 
\end{tikzpicture}
\begin{tikzpicture}[scale = 0.128]
\clip(-2.5,-2.5) rectangle (21.5,11.5);
\draw[lightgray, very thin] (-6,-6) grid (22,12) ;
\draw[dotted] (-4,-2) -- (24,12);
\foreach \x in {-1,0,1,2,3,4,5,6} {
\draw[fill] (\x*12,\x*6) circle (12 pt);
\draw[fill,shift={(2,1)}] (\x*12,\x*6) circle (12 pt);
\draw[shift={(6,3)}] (\x*12,\x*6) circle (12 pt);
\draw[shift={(8,4)}] (\x*12,\x*6) circle (12 pt);}
\end{tikzpicture}
\begin{tikzpicture}[scale = 0.128]
\clip(-2.5,-2.5) rectangle (21.5,11.5);
\draw[lightgray, very thin] (-6,-6) grid (22,12) ;
\draw[dotted] (-4,-2) -- (24,12);
\foreach \x in {-1,0,1,2,3,4,5,6} {
\draw[fill] (\x*6,\x*3) circle (12 pt);
\draw[shift={(2,1)}] (\x*6,\x*3) circle (12 pt);}
\end{tikzpicture}

\vspace*{0.1 cm}
\begin{tikzpicture}[scale = 0.128]
\clip(-2.5,-2.5) rectangle (21.5,11.5);
\draw[lightgray, very thin] (-6,-6) grid (22,12) ;
\draw[dotted] (-4,-2) -- (24,12);
\draw[dotted] (-3,2) -- (17,12);
\foreach \x in {-1,0,1,2,3,4,5,6}{
\draw[fill] (\x*4,\x*2) circle (12 pt);
\draw[shift={(1,4)}] (\x*4,\x*2) circle (12 pt);
};
\end{tikzpicture}
\begin{tikzpicture}[scale = 0.128]
\clip(-2.5,-2.5) rectangle (21.5,11.5);
\draw[lightgray, very thin] (-6,-6) grid (22,12) ;
\draw[dotted] (-4,-2) -- (24,12);
\draw[dotted] (-3,2) -- (17,12);
\foreach \x in {-1,0,1,2,3,4,5,6}{
\draw[fill] (\x*8,\x*4) circle (12 pt);
\draw[shift={(4,2)}] (\x*8,\x*4) circle (12 pt);
\draw[fill,shift={(1,4)}] (\x*8,\x*4) circle (12 pt);
\draw[shift={(5,6)}] (\x*8,\x*4) circle (12 pt);
};
\end{tikzpicture}
\begin{tikzpicture}[scale = 0.128]
\clip(-2.5,-2.5) rectangle (21.5,11.5);
\draw[lightgray, very thin] (-6,-6) grid (22,12) ;
\draw[dotted] (-4,-2) -- (24,12);
\draw[dotted] (-5,0) -- (23,14);
\foreach \x in {-1,0,1,2,3,4} {
\draw[fill] (\x*12,\x*6) circle (12 pt);
\draw[fill,shift={(2,1)}] (\x*12,\x*6) circle (12 pt);
\draw[shift={(6,3)}] (\x*12,\x*6) circle (12 pt);
\draw[shift={(8,4)}] (\x*12,\x*6) circle (12 pt);
\draw[fill,shift={(-1,2)}] (\x*12,\x*6) circle (12 pt);
\draw[fill,shift={(1,3)}] (\x*12,\x*6) circle (12 pt);
\draw[shift={(5,5)}] (\x*12,\x*6) circle (12 pt);
\draw[shift={(7,6)}] (\x*12,\x*6) circle (12 pt);}
\end{tikzpicture}
\begin{tikzpicture}[scale = 0.128]
\clip(-2.5,-2.5) rectangle (21.5,11.5);
\draw[lightgray, very thin] (-6,-6) grid (22,12) ;
\draw[dotted] (-4,-2) -- (24,12);
\draw[dotted] (-5,0) -- (23,14);
\foreach \x in {-1,0,1,2,3,4} {
\draw[fill] (\x*12,\x*6) circle (12 pt);
\draw[fill,shift={(2,1)}] (\x*12,\x*6) circle (12 pt);
\draw[fill,shift={(6,3)}] (\x*12,\x*6) circle (12 pt);
\draw[fill,shift={(8,4)}] (\x*12,\x*6) circle (12 pt);
\draw[shift={(-1,2)}] (\x*12,\x*6) circle (12 pt);
\draw[shift={(1,3)}] (\x*12,\x*6) circle (12 pt);
\draw[shift={(5,5)}] (\x*12,\x*6) circle (12 pt);
\draw[shift={(7,6)}] (\x*12,\x*6) circle (12 pt);}
\end{tikzpicture}

\vspace*{0.1 cm}
\begin{tikzpicture}[scale = 0.128]
\clip(-2.5,-2.5) rectangle (21.5,11.5);
\draw[lightgray, very thin] (-6,-6) grid (22,12) ;
\draw[dotted] (-4,-2) -- (24,12);
\draw[dotted] (-5,0) -- (23,14);
\foreach \x in {-1,0,1,2,3,4} {
\draw[fill] (\x*12,\x*6) circle (12 pt);
\draw[fill,shift={(2,1)}] (\x*12,\x*6) circle (12 pt);
\draw[shift={(6,3)}] (\x*12,\x*6) circle (12 pt);
\draw[shift={(8,4)}] (\x*12,\x*6) circle (12 pt);
\draw[shift={(-1,2)}] (\x*12,\x*6) circle (12 pt);
\draw[shift={(1,3)}] (\x*12,\x*6) circle (12 pt);
\draw[fill,shift={(5,5)}] (\x*12,\x*6) circle (12 pt);
\draw[fill,shift={(7,6)}] (\x*12,\x*6) circle (12 pt);}
\end{tikzpicture}
\begin{tikzpicture}[scale = 0.128]
\clip(-2.5,-2.5) rectangle (21.5,11.5);
\draw[lightgray, very thin] (-6,-6) grid (22,12) ;
\draw[dotted] (-4,-2) -- (24,12);
\draw[dotted] (-5,0) -- (23,14);
\foreach \x in {-1,0,1,2,3,4} {
\draw[fill] (\x*12,\x*6) circle (12 pt);
\draw[shift={(2,1)}] (\x*12,\x*6) circle (12 pt);
\draw[fill,shift={(6,3)}] (\x*12,\x*6) circle (12 pt);
\draw[shift={(8,4)}] (\x*12,\x*6) circle (12 pt);
\draw[fill,shift={(-1,2)}] (\x*12,\x*6) circle (12 pt);
\draw[shift={(1,3)}] (\x*12,\x*6) circle (12 pt);
\draw[fill,shift={(5,5)}] (\x*12,\x*6) circle (12 pt);
\draw[shift={(7,6)}] (\x*12,\x*6) circle (12 pt);}
\end{tikzpicture}
\begin{tikzpicture}[scale = 0.128]
\clip(-2.5,-2.5) rectangle (21.5,11.5);
\draw[lightgray, very thin] (-6,-6) grid (22,12) ;
\draw[dotted] (-4,-2) -- (24,12);
\draw[dotted] (-5,0) -- (23,14);
\foreach \x in {-1,0,1,2,3,4} {
\draw[fill] (\x*12,\x*6) circle (12 pt);
\draw[shift={(2,1)}] (\x*12,\x*6) circle (12 pt);
\draw[fill,shift={(6,3)}] (\x*12,\x*6) circle (12 pt);
\draw[shift={(8,4)}] (\x*12,\x*6) circle (12 pt);
\draw[shift={(-1,2)}] (\x*12,\x*6) circle (12 pt);
\draw[fill,shift={(1,3)}] (\x*12,\x*6) circle (12 pt);
\draw[shift={(5,5)}] (\x*12,\x*6) circle (12 pt);
\draw[fill,shift={(7,6)}] (\x*12,\x*6) circle (12 pt);}
\end{tikzpicture}
\begin{tikzpicture}[scale = 0.128]
\clip(-2.5,-2.5) rectangle (21.5,11.5);
\draw[lightgray, very thin] (-6,-6) grid (22,12) ;
\draw[dotted] (-4,-2) -- (24,12);
\draw[dotted] (-1,-5) -- (3,15);
\foreach \x in {-1,0,1,2,3,4} 
\foreach \y in {-2,-1,0,1,2,3}
{
\draw[fill] (\x*6 +\y,\x*3+5*\y) circle (12 pt);
\draw[shift ={(1,-1)}] (\x*6 +\y,\x*3+5*\y) circle (12 pt);
};
\end{tikzpicture}

\vspace*{0.1 cm}
\begin{tikzpicture}[scale = 0.128]
\clip(-2.5,-2.5) rectangle (21.5,11.5);
\draw[lightgray, very thin] (-6,-6) grid (22,12) ;
\draw[dotted] (-4,-2) -- (24,12);
\draw[dotted] (-1,-5) -- (3,15);
\foreach \x in {-1,0,1,2,3} 
\foreach \y in {-2,-1,0,1,2}
{
\draw[fill] (\x*6 +2*\y,\x*3+10*\y) circle (12 pt);
\draw[fill,shift ={(2,1)}] (\x*6 +2*\y,\x*3+10*\y) circle (12 pt);
\draw[shift ={(1,5)}] (\x*6 +2*\y,\x*3+10*\y) circle (12 pt);
\draw[shift ={(3,6)}] (\x*6 +2*\y,\x*3+10*\y) circle (12 pt);
};
\end{tikzpicture}
\begin{tikzpicture}[scale = 0.128]
\clip(-2.5,-2.5) rectangle (21.5,11.5);
\draw[lightgray, very thin] (-6,-6) grid (22,12) ;
\draw[dotted] (-4,-2) -- (24,12);
\draw[dotted] (-1,-5) -- (3,15);
\foreach \x in {-1,0,1,2,3} 
\foreach \y in {-2,-1,0,1,2}
{
\draw[fill] (\x*6 +2*\y,\x*3+10*\y) circle (12 pt);
\draw[shift ={(2,1)}] (\x*6 +2*\y,\x*3+10*\y) circle (12 pt);
\draw[fill,shift ={(1,5)}] (\x*6 +2*\y,\x*3+10*\y) circle (12 pt);
\draw[shift ={(3,6)}] (\x*6 +2*\y,\x*3+10*\y) circle (12 pt);
};
\end{tikzpicture}
\begin{tikzpicture}[scale = 0.128]
\clip(-2.5,-2.5) rectangle (21.5,11.5);
\draw[lightgray, very thin] (-6,-6) grid (22,12) ;
\draw[dotted] (-4,-2) -- (24,12);
\draw[dotted] (-1,-5) -- (3,15);
\foreach \x in {-1,0,1,2,3} 
\foreach \y in {-2,-1,0,1,2}
{
\draw[fill] (\x*6 +2*\y,\x*3+10*\y) circle (12 pt);
\draw[shift ={(2,1)}] (\x*6 +2*\y,\x*3+10*\y) circle (12 pt);
\draw[shift ={(1,5)}] (\x*6 +2*\y,\x*3+10*\y) circle (12 pt);
\draw[fill,shift ={(3,6)}] (\x*6 +2*\y,\x*3+10*\y) circle (12 pt);
};
\end{tikzpicture}
\begin{tikzpicture}[scale = 0.128]
\clip(-2.5,-1.5) rectangle (21.5,12.5);
\draw[lightgray, very thin] (-6,-6) grid (22,12) ;
\draw[dotted] (-4,-2) -- (24,12);
\draw[dotted] (1,-2) -- (-6,12);
\foreach \x in {0,1} 
\foreach \y in {-1,0}
{
\draw[fill] (\x*12- \y*6,\x*6+\y*12) circle (12 pt);
\draw[fill,shift={(2,1)}] (\x*12- \y*6,\x*6+\y*12) circle (12 pt);
\draw[fill,shift={(6,3)}] (\x*12- \y*6,\x*6+\y*12) circle (12 pt);
\draw[fill,shift={(8,4)}] (\x*12- \y*6,\x*6+\y*12) circle (12 pt);
\draw[fill,shift={(-1,2)}] (\x*12- \y*6,\x*6+\y*12) circle (12 pt);
\draw[fill,shift={(1,3)}] (\x*12- \y*6,\x*6+\y*12) circle (12 pt);
\draw[fill,shift={(5,5)}] (\x*12- \y*6,\x*6+\y*12) circle (12 pt);
\draw[fill,shift={(7,6)}] (\x*12- \y*6,\x*6+\y*12) circle (12 pt);
\draw[shift={(-3,6)}] (\x*12- \y*6,\x*6+\y*12) circle (12 pt);
\draw[shift={(-1,7)}] (\x*12- \y*6,\x*6+\y*12) circle (12 pt);
\draw[shift={(3,9)}] (\x*12- \y*6,\x*6+\y*12) circle (12 pt);
\draw[shift={(5,10)}] (\x*12- \y*6,\x*6+\y*12) circle (12 pt);
\draw[shift={(-4,8)}] (\x*12- \y*6,\x*6+\y*12) circle (12 pt);
\draw[shift={(-2,9)}] (\x*12- \y*6,\x*6+\y*12) circle (12 pt);
\draw[shift={(2,11)}] (\x*12- \y*6,\x*6+\y*12) circle (12 pt);
\draw[shift={(4,12)}] (\x*12- \y*6,\x*6+\y*12) circle (12 pt);
};
\end{tikzpicture}

\vspace*{0.1 cm}
\begin{tikzpicture}[scale = 0.128]
\clip(-2.5,-1.5) rectangle (21.5,12.5);
\draw[lightgray, very thin] (-6,-6) grid (22,12) ;
\draw[dotted] (-4,-2) -- (24,12);
\draw[dotted] (1,-2) -- (-6,12);
\foreach \x in {0,1} 
\foreach \y in {-1,0}
{
\draw[fill] (\x*12- \y*6,\x*6+\y*12) circle (12 pt);
\draw[fill,shift={(2,1)}] (\x*12- \y*6,\x*6+\y*12) circle (12 pt);
\draw[shift={(6,3)}] (\x*12- \y*6,\x*6+\y*12) circle (12 pt);
\draw[shift={(8,4)}] (\x*12- \y*6,\x*6+\y*12) circle (12 pt);
\draw[fill,shift={(-1,2)}] (\x*12- \y*6,\x*6+\y*12) circle (12 pt);
\draw[fill,shift={(1,3)}] (\x*12- \y*6,\x*6+\y*12) circle (12 pt);
\draw[shift={(5,5)}] (\x*12- \y*6,\x*6+\y*12) circle (12 pt);
\draw[shift={(7,6)}] (\x*12- \y*6,\x*6+\y*12) circle (12 pt);
\draw[shift={(-3,6)}] (\x*12- \y*6,\x*6+\y*12) circle (12 pt);
\draw[shift={(-1,7)}] (\x*12- \y*6,\x*6+\y*12) circle (12 pt);
\draw[fill,shift={(3,9)}] (\x*12- \y*6,\x*6+\y*12) circle (12 pt);
\draw[fill,shift={(5,10)}] (\x*12- \y*6,\x*6+\y*12) circle (12 pt);
\draw[shift={(-4,8)}] (\x*12- \y*6,\x*6+\y*12) circle (12 pt);
\draw[shift={(-2,9)}] (\x*12- \y*6,\x*6+\y*12) circle (12 pt);
\draw[fill,shift={(2,11)}] (\x*12- \y*6,\x*6+\y*12) circle (12 pt);
\draw[fill,shift={(4,12)}] (\x*12- \y*6,\x*6+\y*12) circle (12 pt);
};
\end{tikzpicture}
\begin{tikzpicture}[scale = 0.128]
\clip(-2.5,-1.5) rectangle (21.5,12.5);
\draw[lightgray, very thin] (-6,-6) grid (22,12) ;
\draw[dotted] (-4,-2) -- (24,12);
\draw[dotted] (1,-2) -- (-6,12);
\foreach \x in {0,1} 
\foreach \y in {-1,0}
{
\draw[fill] (\x*12- \y*6,\x*6+\y*12) circle (12 pt);
\draw[fill,shift={(2,1)}] (\x*12- \y*6,\x*6+\y*12) circle (12 pt);
\draw[fill,shift={(6,3)}] (\x*12- \y*6,\x*6+\y*12) circle (12 pt);
\draw[fill,shift={(8,4)}] (\x*12- \y*6,\x*6+\y*12) circle (12 pt);
\draw[shift={(-1,2)}] (\x*12- \y*6,\x*6+\y*12) circle (12 pt);
\draw[shift={(1,3)}] (\x*12- \y*6,\x*6+\y*12) circle (12 pt);
\draw[shift={(5,5)}] (\x*12- \y*6,\x*6+\y*12) circle (12 pt);
\draw[shift={(7,6)}] (\x*12- \y*6,\x*6+\y*12) circle (12 pt);
\draw[fill,shift={(-3,6)}] (\x*12- \y*6,\x*6+\y*12) circle (12 pt);
\draw[fill,shift={(-1,7)}] (\x*12- \y*6,\x*6+\y*12) circle (12 pt);
\draw[fill,shift={(3,9)}] (\x*12- \y*6,\x*6+\y*12) circle (12 pt);
\draw[fill,shift={(5,10)}] (\x*12- \y*6,\x*6+\y*12) circle (12 pt);
\draw[shift={(-4,8)}] (\x*12- \y*6,\x*6+\y*12) circle (12 pt);
\draw[shift={(-2,9)}] (\x*12- \y*6,\x*6+\y*12) circle (12 pt);
\draw[shift={(2,11)}] (\x*12- \y*6,\x*6+\y*12) circle (12 pt);
\draw[shift={(4,12)}] (\x*12- \y*6,\x*6+\y*12) circle (12 pt);
};
\end{tikzpicture}
\begin{tikzpicture}[scale = 0.128]
\clip(-2.5,-1.5) rectangle (21.5,12.5);
\draw[lightgray, very thin] (-6,-6) grid (22,12) ;
\draw[dotted] (-4,-2) -- (24,12);
\draw[dotted] (1,-2) -- (-6,12);
\foreach \x in {0,1} 
\foreach \y in {-1,0}
{
\draw[fill] (\x*12- \y*6,\x*6+\y*12) circle (12 pt);
\draw[fill,shift={(2,1)}] (\x*12- \y*6,\x*6+\y*12) circle (12 pt);
\draw[shift={(6,3)}] (\x*12- \y*6,\x*6+\y*12) circle (12 pt);
\draw[shift={(8,4)}] (\x*12- \y*6,\x*6+\y*12) circle (12 pt);
\draw[shift={(-1,2)}] (\x*12- \y*6,\x*6+\y*12) circle (12 pt);
\draw[shift={(1,3)}] (\x*12- \y*6,\x*6+\y*12) circle (12 pt);
\draw[fill,shift={(5,5)}] (\x*12- \y*6,\x*6+\y*12) circle (12 pt);
\draw[fill,shift={(7,6)}] (\x*12- \y*6,\x*6+\y*12) circle (12 pt);
\draw[fill,shift={(-3,6)}] (\x*12- \y*6,\x*6+\y*12) circle (12 pt);
\draw[fill,shift={(-1,7)}] (\x*12- \y*6,\x*6+\y*12) circle (12 pt);
\draw[shift={(3,9)}] (\x*12- \y*6,\x*6+\y*12) circle (12 pt);
\draw[shift={(5,10)}] (\x*12- \y*6,\x*6+\y*12) circle (12 pt);
\draw[shift={(-4,8)}] (\x*12- \y*6,\x*6+\y*12) circle (12 pt);
\draw[shift={(-2,9)}] (\x*12- \y*6,\x*6+\y*12) circle (12 pt);
\draw[fill,shift={(2,11)}] (\x*12- \y*6,\x*6+\y*12) circle (12 pt);
\draw[fill,shift={(4,12)}] (\x*12- \y*6,\x*6+\y*12) circle (12 pt);
};
\end{tikzpicture}
\begin{tikzpicture}[scale = 0.128]
\clip(-2.5,-1.5) rectangle (21.5,12.5);
\draw[lightgray, very thin] (-6,-6) grid (22,12) ;
\draw[dotted] (-4,-2) -- (24,12);
\draw[dotted] (1,-2) -- (-6,12);
\foreach \x in {0,1} 
\foreach \y in {-1,0}
{
\draw[fill] (\x*12- \y*6,\x*6+\y*12) circle (12 pt);
\draw[shift={(2,1)}] (\x*12- \y*6,\x*6+\y*12) circle (12 pt);
\draw[fill,shift={(6,3)}] (\x*12- \y*6,\x*6+\y*12) circle (12 pt);
\draw[shift={(8,4)}] (\x*12- \y*6,\x*6+\y*12) circle (12 pt);
\draw[shift={(-1,2)}] (\x*12- \y*6,\x*6+\y*12) circle (12 pt);
\draw[fill,shift={(1,3)}] (\x*12- \y*6,\x*6+\y*12) circle (12 pt);
\draw[shift={(5,5)}] (\x*12- \y*6,\x*6+\y*12) circle (12 pt);
\draw[fill,shift={(7,6)}] (\x*12- \y*6,\x*6+\y*12) circle (12 pt);
\draw[fill,shift={(-3,6)}] (\x*12- \y*6,\x*6+\y*12) circle (12 pt);
\draw[shift={(-1,7)}] (\x*12- \y*6,\x*6+\y*12) circle (12 pt);
\draw[fill,shift={(3,9)}] (\x*12- \y*6,\x*6+\y*12) circle (12 pt);
\draw[shift={(5,10)}] (\x*12- \y*6,\x*6+\y*12) circle (12 pt);
\draw[shift={(-4,8)}] (\x*12- \y*6,\x*6+\y*12) circle (12 pt);
\draw[fill,shift={(-2,9)}] (\x*12- \y*6,\x*6+\y*12) circle (12 pt);
\draw[shift={(2,11)}] (\x*12- \y*6,\x*6+\y*12) circle (12 pt);
\draw[fill,shift={(4,12)}] (\x*12- \y*6,\x*6+\y*12) circle (12 pt);
};
\end{tikzpicture}
\caption{20 classes of possible choices of $V_+,V_-$ in 
Corollary \ref{Cor:Z2}, where every class is illustrated by a typical representative. 
Here $\Z^2$ is embedded into $\R^2$ in the canonical way. 
The vertices of $V_+$ are marked in black, 
the vertices of $V_-$ are marked in white. 
The origin $o$ is in one of the vertices of $V_+$.  
For classes 3-20 the pattern extends periodically. 
For some of the classes the set $V_\pm = V_+ \cup V_-$ is the same: for classes 1-2 the set consists of 4 vertices, for 3-4 all vertices are on a single line, for 5-6  and for 7-11 all vertices are on two parallel lines, 
 for 12,13-15 and 16-20 the vertices are on infinitely many parallel lines. Case (a) of the corollary corresponds to the classes 4,5, and case (b) to 12,14. 
}
\label{Fig:Z2}
\end{figure}
In fact there are many more classes of choices of suitable $V_+,V_-$ satisfying our symmetry assumptions; some of them are illustrated in Figure \ref{Fig:Z2}. But even simple classs  may be of interest: 
E.g. the choice (a) in Corollary \ref{Cor:Z2}
implies that for subcritical Bernoulli percolation 
on $\Z^2$, the expected number of vertices 
of the cluster of $(0,0)$ in the layer $V_n = \{(n,k) : k \in \Z\}$ 
is maximal for $n = 0$. This is a very intuitive assertion 
on the shape of clusters, which seems to be novel and may be of 
interest. 
%
%

Our last application concerns hypercubes. 
Here we formulate the conclusion concerning the expectations in terms of 
connection probabilities. 

\begin{Cor} \label{Cor:hypercube} Hypercube graphs. 
For $d \ge 1$ we consider the $d$-dimensional hypercube graph 
$L_2^d = L_2 \cp ... \cp L_2$ and $o = (0,...,0)$. 
Let $c_i(p) := \P_p(o \lra v)$ where $d(v,o) = i$. 
For all $k,l \ge 0$ such that $k+l \le d$ we have 
$$
\sum_{i=0}^k \sum_{j=0}^l c_{i+j} (-1)^i \binom k i \binom l j \ge 0
\; \text{ and } \; 
| \sum_{i=0}^k c_{i} (-1)^i \binom k i | \ge 
| \sum_{i=0}^k c_{i+l} (-1)^i \binom k i |.
$$
In particular for the discrete derivatives of the sequence $c = (c_0,...,c_d)$  we have 
$$
\forall k \le d: (-1)^k \De^k[c](0) \ge 0 \quad \text{ and }  \quad 
\forall l \le d-k: |\De^k[c](0)| \ge |\De^k[c](l)|.
$$
(Here: $\De [a_0,...,a_m] = [a_1-a_0,...,a_m - a_{m-1}]$ and $\De^{k+1}[a] = \De[\De^k[a]]$.)
\end{Cor}
Similarly other relations between the connection probabilities can be obtained for hypercube graphs, 
but the above are the most interesting to us. 
We note that the validity of the bunkbed conjecture for 
hypercube graphs would imply that 
$c_i \le c_{i-1}$ for all $i$, i.e. $\De^1[c] (i) \le 0$
for all $i$. While the bunkbed conjecture thus concerns the first derivatives of the sequence at all points, the above corollary 
concerns all derivatives of the sequence at $0$. 
We note that similarly to the above corollary we can obtain 
relations between connection probabilities in other highly symmetric graphs such as the graphs of (quasi-) regular polyhedra, 
symmetric complete bipartite graphs or the Petersen graph and its relatives.

\section{Group actions on pairs of sets} \label{Sec:action}

A major ingredient of the proof of our main theorem 
is a combinatorial result concerning group actions on pairs of sets, 
which will be presented in this section.  
For some basics and results on group actions in the context of graph theory see \cite{GR}. 
We start by reviewing properties of group actions and corresponding notations. 
Let $\G$ be a group acting on a set $X$. 
Our notation for the group action will be $g(x)$ for $g \in \G$ and $x \in X$. 
The action of $\G$ on $X$ naturally extends to subsets
via $g(A) := \{g(x): x \in A\}$ for $g \in \G, A \subset X$.  
The action of $\G$ on two sets $X,X'$ naturally extends to 
pairs via $g(x,x') = (g(x),g(x'))$ for $g \in \G, x \in X,x' \in X'$. 
We will use standard notation for orbits and stabilizers: 
For a subset $\G' \subset \G$ and for $x \in X$ 
we write $\G'(x) = \{g(x): g \in \G'\}$ 
and similar for sets and pairs of sets. 
Also we write $\G_x =  \{g \in \G: g(x) = x\}$ for given 
$x \in X$ and similar for sets. 
Here the use of multiple indices means that all elements/sets have to be stabilized, 
e.g. $\G_{x,A} = \{g \in \G: g(x) = x, g(A) = A\}$. 
We say that $\G$ acts transitively on $X$ iff 
$\forall x,y \in X \exists g \in \G: g(x) = y$. 

\begin{Thm} \label{Prop:combi}
Let $X,X'$ be sets, and let the set of pairs of finite subsets of $X$ and $X'$ be denoted by 
$\bar M = \{(A,A'): A \subset X,A' \subset X', |A|,|A'| < \infty\}$. 
Let $\G$ be a group acting transitively 
both on $X$ and on $X'$, and suppose that 
\begin{align} \label{equ:stabsym}
\forall x \in X,x' \in X': |\G_x(x')| = |\G_{x'}(x)|.
\end{align}
Then for every $\bar A  = (A,A') \in \bar M$ and all 
$o \in X$, $o' \in X'$ 
\begin{align*}
&\text{ for } \quad 
[\bar{A}]_{o} := \{\bar B \in  \G(\bar A): o \in B\} 
\quad \text{ and } \quad [\bar{A}]_{o'} := \{\bar B \in \G(\bar A), o' \in B'\} \\
&\text{ we have } \quad |[\bar{A}]_{o}| \cdot |A'| = |[\bar{A}]_{o'}| \cdot|A|. 
\end{align*}
If additionally the orbits in \eqref{equ:stabsym}
are finite, then also $|[\bar{A}]_{o}|,|[\bar{A}]_{o'}|  < \infty$.
\end{Thm}

\Pf Let us first assume that $|A|,|A'| > 0$.
We observe that 
$$
[\bar A]_{o} = \bigcup_{x \in A} [\bar A]^{x}_{o} 
\quad \text{ for } [\bar A]^{x}_{o} := \{g(\bar A): g \in \G, g(x) = o\}.  
$$
This union is not disjoint. 
However we have $[\bar A]^{x}_{o} = [\bar A]^{y}_{o}$ for $y \in \G_{\bar A}(x)$, 
and $[\bar A]^{x}_{o} \cap [\bar A]^{y}_{o} = \emptyset$ otherwise. Noting that $\G_{\bar A}(x) \subset A$ is finite,
we thus obtain  
$$
|[\bar A]_o| = \sum_{x \in A} \frac{|[\bar A]^{x}_{o}|}{|\G_{\bar A}(x)|}.
$$
Next for every $x \in A$ we choose $g_{x} \in \G$ 
such that $g_{x}(x) = o$ (using the transitivity) 
and note that $\{g \in \G: g(x) = o\} = g_x \G_x$. 
This gives $[\bar A]^{x}_{o} = g_{x}(\G_{x}(\bar A))$, 
which implies $|[\bar A]^{x}_{o}| = |\G_{x}(\bar A)|$. 
In order to obtain a suitable expression for the last term, 
we fix some $x' \in A'$. 
The following lemma (used twice) implies 
$$
|\G_{x}(\bar A)| \cdot |\G_{x,\bar A}(x')| 
= |\G_{x}((\bar A,x'))| =  
|\G_{x}(x')| \cdot |\G_{x,x'}(\bar A)|.
$$ 
Noting that $\G_{x,\bar A}(x') \subset A'$ is  also finite, 
we may combine the above equalities we get  
$$
|[\bar A]_o| = \sum_{x \in A} \frac{|\G_x(\bar A)|}{|\G_{\bar A}(x)|}
= \sum_{x \in A} \frac{|\G_x(x')| \cdot |\G_{x,x'}(\bar A)|}{|\G_{\bar A}(x)| \cdot |\G_{x,\bar A}(x')|}, 
$$
Again by the following lemma we have 
$|\G_{\bar A}(x)| \cdot |\G_{x,\bar A}(x')| = |\G_{\bar A}((x,x'))|$. 
Summing both sides over all possible choices of $x' \in A'$
we obtain 
\begin{align*}
&|[\bar A]_o| \cdot |A'| = 
\sum_{x \in A,x' \in A'} \frac{|\G_x(x')| \cdot |\G_{x,x'}(\bar A)|}{|\G_{\bar A}((x,x'))|}, \quad \text{ and similarly we have }\\
&|[\bar A]_{o'}| \cdot |A| =  \sum_{x \in A,x' \in A'} \frac{|\G_{x'}(x)| \cdot |\G_{x,x'}(\bar A)|}{|\G_{\bar A}((x,x'))|}. 
\end{align*}
Since we have assumed that $|\G_{x'}(x)| = |\G_{x}(x')|$, 
the above sums are equal,
which gives the first part of the claim. 
If the orbits of stabilizers are finite, 
then also  
\begin{align*}
|\G_{x,x'}(\bar A)| &\le 
|\G_{x,x'}(A)| \cdot |\G_{x,x'}(A')| \le 
|\G_{x'}(A)| \cdot |\G_{x}(A')| \\
&\le 
\prod_{y \in A} |\G_{x'}(y)|
\prod_{y' \in A'}
|\G_{x}(y')|  < \infty,  
\end{align*}
and thus the above sums are finite, which gives the second part of the claim. 

\newpage 

Finally, in case of $|A| = 0$ or $|A'| = 0$ we note that
$|A| = 0$ implies that $[\bar A]_o = 0$ 
(and similarly for $A'$), so in this case the first claim 
holds trivially. 
We also note that 
$|[(A,\emptyset)]_o| \le |[(A,A')]_o|$ for all 
finite $A \subset X, A' \subset X'$ (and similarly for $A'$), 
so in this case the second claim 
follows from the finiteness in the main case. 
\qed

\begin{Lem} \label{Lem:orbits}
Let $\G$ be a group acting on the sets $X$ and $X'$. 
For $x \in X$, $x' \in X'$
we have  
$|\G((x,x'))| = |\G(x)| \cdot |\G_x(x')|.$
\end{Lem}

\Pf Let $y \in \G(x)$, i.e. $y = g_0(x)$ for some $g_0 \in \G$.
We have
\begin{align*}
&\{g(x'): g \in \G \text{ s.t. } g(x) = y\} 
= \{g(x'): g \in g_0 \G_x\} 
= g_0(\G_x(x')), 
\; \text{  which implies}\\
&|\G((x,x'))| = \sum_{y \in \G(x)} |\{g(x'): g \in \G \text{ s.t. } g(x) = y\}|= |\G(x)| \cdot |\G_x(x')|.
\end{align*}

\vspace*{ - 1.2 cm}
\qed 

\smallskip 

\begin{Rem} Concerning condition \eqref{equ:stabsym}. 
\begin{itemize}
\item 
In case of $\bar{A} = (A,A') = (\{x\},\{x'\})$ for $x \in X, x' \in X'$ we note that $|[\bar A]_o| = |\{g(x') : g \in \G \text{ s.t. }g(x) = o\}| = |\G_x(x')|$ since the action is transitive on $X$, 
and similarly   $|[\bar A]_{o'}|  = |\G_{x'}(x)|$. 
So in this special case the claim is equivalent to $|\G_x(x')| = |\G_{x'}(x)|$. In that sense condition \eqref{equ:stabsym} can be seen to be very natural. 
\item  
Condition \eqref{equ:stabsym} is not very explicit. 
For more explicit stronger conditions see the next section. 
Under the stronger assumption  $|X| = |X'| < \infty$ 
the theorem can be shown with a less complicated double 
counting argument. 
\end{itemize} 
\end{Rem}


\section{Random partitions} \label{Sec:partitions}

Here we formulate and prove a more general version of our main result for random partitions of a set. 
We first introduce the setting. Let $X$ be a countable set. 
Let $P_X$ denote the set of all partitions of $X$ into disjoint subsets. 
For $x,y \in X$ let $\{x \lra y\}$ denote the set of 
all partitions such that $x,y$ belong to the same set of the partition. 
$\F_X := \si(\{x \lra y\}: x,y \in X)$ is a suitable $\si$-algebra on $P_X$. 
A probability distribution $\P$ on $(P_X,\F_X)$ will 
be called a random partition of $X$. 
For $x \in X$ and $P \in P_X$ let  
$\cC_x(P)$ denote the set of the partition $P$ containing $x$, 
so that $y \in \cC_x(P)$ iff $P \in \{x \lra y\}$. 
$\cC_x$ is a random variable on $(P_X,\F_X)$, called the cluster of $x$. 
Suitable symmetry assumptions will be formulated 
via the action of a group $\G$ on $X$. 
This action naturally extends  to $P_X$ via 
$g(P) = \{g(A): A \in P\}$ for $P \in P_X$. 
We note that every $g \in \G$ may thus be considered 
a measurable map 
from $(P_X,\F_X)$ to $(P_X,\F_X)$.

\begin{Thm} \label{Thm:partitions} 
Let $X_+,X_-$ be disjoint countable sets and 
let $\G$ be a group acting on $X := X_+ \cup X_-$ 
such that 
\begin{itemize}
\item[$(\G 1')$] $\forall g \in \G: g(X_+),g(X_-) \in \{X_+,X_-\}$,
\item[$(\G 2')$] 
$\G$ acts transitively on $X$ and 
\item[$(\G 3')$]
$\forall x_+ \in X_+,x_- \in X_-: \quad 
|\G_{x_+}(x_-)| = |\G_{x_-}(x_+)|< \infty$.
\end{itemize}
Let $\P$ be a random partition of $X$, let $\cC := \cC_o$ denote the cluster of some fixed $o \in X_+$ and let  $\cC_+ := \cC \cap X_+$ and $\cC_- := \cC \cap X_-$.
Suppose that $\P$ is invariant under every $g \in \G$.  
Then for every bounded $f: \N  \to \R$ we have 
\begin{align*}
\E\Big(\!\big(f(|\cC_+|) \!-\! f(|\cC_-|)\big)1_{\{|\cC| < \infty\}}\!\Big) = \E\Big(\frac{(|\cC_+| \!-\! |\cC_-|)(f(|\cC_+|) \!-\! f(|\cC_-|))}{|\cC_+| + |\cC_-|}1_{\{|\cC| < \infty\}}\!\Big). 
\end{align*}
In particular we have $|\cC_-| \preceq |\cC_+|$ w.r.t. $\P(.||\cC| < \infty)$ (if $\P(|\cC| < \infty) >0$).  
\end{Thm}

\Pf We have $X_+ \neq \emptyset$ and w.l.o.g. also $X_- \neq \emptyset$. 
By ($\G 1'$)
\begin{align*}
\G = \G_+ \cup \G_-, \quad \text{ where } \quad 
&\G_+ = \{g \in \G: g(X_+) = X_+,g(X_-) = X_-\} \\
\text{ and } \quad 
&\G_- = \{g \in \G: g(X_+) = X_-,g(X_-) = X_+\}.  
\end{align*}
($\G 2'$)  implies $\G_- \neq \emptyset$, 
and fixing some $\tau \in \G_-$ we have $\G_- = \tau \G_+$. 
Let 
\begin{align*}
&M := \{\bar A \subset X:  |\bar A| < \infty\},  \quad 
M_o := \{\bar A \in M: o \in A\} 
\quad \text{ and for } \bar A  \in M_o \text{ let}\\
&[\bar A]_+ := \G_+(\bar A) \cap M_o, \quad  
[\bar A]_- := \G_-(\bar A) \cap M_o \quad \text{ and } \quad 
[\bar A] := \G(\bar A) \cap M_o,  
\end{align*}
so that $[\bar A] = [\bar A]_+ \cup [\bar A]_-.$
We now would like to apply Theorem \ref{Prop:combi}
to the group action $\G_+$ on the sets $X_+$ and $X_-$. 
By ($\G 2'$) both the action on $X_+$ and the action on $X_-$ is transitive, and since 
$\G_x = (\G_+)_x$ for all $x \in X$, condition \eqref{equ:stabsym} is an immediate consequence of ($\G 3'$). 
Identifying $\bar A = (A_+,A_-) \in \bar M$ with 
$\bar A = A_+ \cup A_- \in M$ 
(where $A_+ \subset X_+,A_- \subset X_-$), we note that 
\begin{align*}
&[\bar A]_+ = \{B_+ \cup B_-: (B_+,B_-)  \in [\bar A]_o\}
\quad \text{ and }\\
&[\bar A]_- = \tau \G_+(\bar A) \cap M_o = 
\{\tau(B_+ \cup B_-): (B_+,B_-) \in [\bar  A]_{\tau^{-1}(o)}\}. 
\end{align*} 
In particular $|[\bar A]_+| = |[\bar A]_o|$ and  $|[\bar A]_-| =| [\bar A]_{\tau^{-1}(o)}|$, thus 
Theorem \ref{Prop:combi} gives 
$$
|[\bar A]_+| \cdot |A_-| = |[\bar A]_-| \cdot |A_+|
\quad \text{ and } \quad 
|[\bar A]_+|, |[\bar A]_-| < \infty \quad \text{ for every $A \in M_o$.}
$$
We now observe that 
for bounded $f: \N  \to \R$ we have 
\begin{align*}
&
\E\big((f(|\cC_+|) - f(|\cC_-|))1_{\{|\cC| < \infty\}}\big)
= \sum_{\bar B \in M_o}  (f(|B_+|)-f(|B_-|)) \P(\cC = \bar{B})\\
&= \sum_{\bar A, \bar B \in M_o: \bar A \in [\bar B]}  \!
\frac {f(|B_+|)-f(|B_-|)} {|[\bar B]|}  \P(\cC = \bar{B}) \\
&= \sum_{\bar A \in M_o} \P(\cC = \bar{A}) \! 
 \sum_{\bar B \in [\bar A]}\frac {f(|B_+|)-f(|B_-|)} {|[\bar A]|}
.
\end{align*}
Here in the first step we have used that on $|\cC| < \infty$ 
we have $\cC \in M_o$ and that $M_o$ is countable, 
since $X$ is countable; 
in the second step we have used that 
$[\bar B]$ is finite for all $\bar B \in M_o$ by the above;  
in the third step we have used that $\P$ is 
$\G$-invariant, and that for all $\bar A,\bar B \in M_o$ 
$$
\bar A \in [\bar B] \quad \Lra \quad  
[\bar A] = [\bar B] \quad \Lra \quad  
\bar B \in [\bar A], 
$$  
which follows from the observation that non-disjoint orbits 
have to coincide. 
We have also used Fubini's theorem, which is justified since 
$f$ is bounded. 
We continue to simplify the inner sum:
We note that  for all $\bar B \in [\bar A]_+$ we have 
$|B_+|  = |A_+|$ and $|B_-|  = |A_-|$, whereas for 
 $\bar B \in [\bar A]_-$ we have 
$|B_+|  = |A_-|$ and $|B_-|  = |A_+|$.
In particular for $|A_+| = |A_-|$ the inner sum equals $0$. 
In case of $|A_+| \neq  |A_-|$ the sets $[\bar A]_+,[\bar A]_-$
are disjoint, and we may rewrite the inner sum 
\begin{align*}
&\sum_{\bar{B} \in [\bar{A}]}\!  \frac{f(|B_+|)-f(|B_-|)}{|[\bar A]|} 
= \sum_{\bar{B} \in [A]_+} \!\! \frac{f(|A_+|)-f(|A_-|)}{|[\bar A]|}
+ \! \sum_{\bar{B} \in [A]_-} \!\! \frac{f(|A_-|)-f(|A_+|)}{|[\bar A]|} \\
&= (f(|A_+|)-f(|A_-|)) \frac{|[\bar A]_+| - |[\bar A]_-|}{|[\bar A]_+| + |[\bar A]_-|}
= (f(|A_+|)-f(|A_-|)) \frac{|A_+|-|A_-|}{|A_+|+|A_-|}. 
\end{align*}
In the last step we have used $\frac{|[\bar A]_-|}{|[\bar A]_+|} = \frac{|A_-|}{|A_+|}$ from above. We have thus shown  
\begin{align*}
&\E\Big(\big(f(|\cC_+|) - f(|\cC_-|)\big)1_{\{|\cC| < \infty\}} )\\ 
&= \sum_{\bar A \in M_o}
\frac{(|A_+|-|A_-|)(f(|A_+|)-f(|A_-|))}{|A_+|+|A_-|} \P_p(\cC = \bar{A}), 
\end{align*}
which equals the claimed expectation. If $f$ is increasing, 
the sum is $\ge 0$, which gives $\E(f(|\cC_+|)1_{\{|\cC| < \infty\}}) \ge \E(f(|\cC_-|)1_{\{|\cC| < \infty\}})$ and thus the claimed domination. \qed

\begin{Rem} 
Some remarks on the above theorem on random partitions:  
\begin{itemize}
\item 
The conditions $(\G 1')$ -  $(\G 3')$ on the action of $\G$ 
should be considered a symmetry assumption on $X_+,X_-$. 
One should think of $X_-$ as a copy of $X_+$ 
with both sets exhibiting a high degree of symmetry. 
\item 
We note that we can obtain similar results in the same way. 
E.g. the same proof shows that under the given conditions we have 
\begin{align*}
\E\Big(\frac{|\cC_-|}{|\cC_+|}1_{\{|\cC| < \infty\}}\Big) = 
\P(|\cC| < \infty, |\cC_-| > 0). 
\end{align*}
\item 
$\P$ is only assumed to be 
$\G$-invariant. 
In particular, one might consider a stochastic process on 
a graph such that the process is invariant under a certain 
subgroup of graph automorphisms, and such that 
the process induces a random partition of the vertex set of the graph. 
For example one could consider the clusters 
of bond percolation with varying edge probabilities 
(say  $p,p' \in (0,1)$ for horizontal and vertical edges 
in $\Z^2$), site percolation, see \cite{G1}, 
or the random cluster model, see \cite{G2}.
\end{itemize}
\end{Rem}
Assumption ($\G 3'$) is not very explicit; in the following lemma we give two natural conditions which imply 
the equality in ($\G 3'$). 
%
%
\begin{Lem} \label{Lem:orbstab}
Let $\G$ be a group acting on $X$ and $x,y \in X$. 
\begin{itemize}
\item[(a)] If for some  $g \in \G$ we have $g(x) = y, g(y) = x$, 
then $|\G_{y}(x)| = |\G_{x}(y)|$.
\item[(b)] If the action is transitive 
and $|X| < \infty$, then  $|\G_{y}(x)| = |\G_{x}(y)|$. 
\end{itemize}\end{Lem} 

\Pf In (a) the claim follows from $\G_{y}(x) = (g^{-1}\G_{x}g)(x) =  g^{-1}(\G_{x}(y))$. 
For (b) we use transitivity and apply  Lemma \ref{Lem:orbits} twice to obtain 
$$
|X| \cdot |\G_{x}(y)|  =| \G((x,y))| =| \G((y,x))| = |X| \cdot |\G_{y}(x)|. 
$$
Together with $|X| < \infty$ this implies  $|\G_{x}(y)| = |\G_{y}(x)|$.  \qed

\section{Independent bond percolation} \label{Sec:percolation}

We now return to the setting of independent percolation from 
the first section and prove the results stated there. 
In order to be able to infer Theorem \ref{Thm:Graphs} from 
Theorem \ref{Thm:partitions}, we we need some auxiliary results.

\begin{Lem} \label{Lem:infinite}
We consider independent bond percolation with parameter $p \in (0,1)$ on a connected, locally finite graph $G = (V,E)$. 
Let $\cC_o$ denote the cluster of a fixed vertex 
$o \in V$.  
Let $V_+,V_- \subset V$ be two disjoint sets of vertices 
such that for some $\tau \in \Aut(G)$ we have $\tau(V_+) = V_-$ and $\tau(V_-) = V+$. 
Then we have $\E_p(|\cC_o \cap V_+|) = \infty$ 
iff $\E_p(|\cC_o \cap V_-|) = \infty$. 
\end{Lem}

\Pf Suppose that $\E_p(|\cC_o \cap V_+|) = \infty$. 
For $o' := \tau(o)$ we have 
\begin{align*}
&\E_p(|\cC_o \cap V_-|) = \E_p(|\cC_{o'} \cap V_+|)
\ge \E_p(|\cC_o \cap V_+| 1_{\{o \lra o'\}}) \\
&= \sum_{n \ge 1} \P_p(|\cC_o \cap V_+| \ge n, o \lra o') 
\ge \sum_{n \ge 1} \P_p(|\cC_o \cap V_+| \ge n)\P_p(o \lra o')\\ 
&=\E_p(|\cC_o \cap V_+|) \P_p(o \lra o') = \infty. 
\end{align*}
In the first step we have used the $\tau$-invariance of $\P_p$, 
in the second step we have used that 
$\cC_{o'} = \cC_o$ on $\{o \lra o'\}$, 
in the third and fifth step we used that 
$\E(X) = \sum_{n \ge 1} \P(X \ge n)$ for any random variable 
with values in $\{0,1,2,...\} \cup \{\infty\}$, 
in the fourth step we have used the FKG-inequality 
(e.g. see \cite{G1}) on the increasing events 
$\{|\cC_o \cap V_+| \ge n\}$ and $\{o \lra o'\}$, 
and the last step follows from 
$\E_p(|\cC_o \cap V_+|) = \infty$ and $\P_p(o \lra o') > 0$.
(We have not used the FKG inequality directly on the expectation, 
since it is usually formulated for random variables that are bounded or have second moments.) 
We have thus shown that $\E_p(|\cC_o \cap V_+|) = \infty$ implies $\E_p(|\cC_o \cap V_-|) = \infty$, and the other implication
is obtained similarly. \qed

\begin{Lem}\label{Lem:orbitstabfinite}
Let $G = (V,E)$ be a connected locally finite graph, 
then $V$ is countable, and for all 
$\G \subset \Aut(G)$ and $x,y \in V$ the orbit $\G_{x}(y)$ is finite. 
\end{Lem}

\Pf For $x \in V$ and $k \ge 0$ let $S_k(x) :=  \{y \in V: d(x,y) = k\}$. We have $|S_0(x)| = 1$ and 
$|S_{k+1}(x)| \le |S_{k}(x)| \cdot \max\{\deg(v): v \in S_k(x)\}$. Since $G$ is locally finite this inductively implies that 
$|S_k(x)| < \infty$ for all $k$. 
In combination with the connectedness of $G$ this shows 
that $V = \bigcup_{k \ge 0} S_k(x)$ is countable
and that for all $x,y \in V$ we have 
$|\G_x(y)| \le |S_{d(x,y)}(x)| < \infty$, noting that graph automorphisms preserve the graph distance. 
\qed 

\bigskip 

{\bf Proof of Theorem \ref{Thm:Graphs}:} 
By Lemma \ref{Lem:orbitstabfinite} $V$ is countable. 
We have $V_+ \neq \emptyset$ and w.l.o.g. 
we may assume $V_- \neq \emptyset$. 
$\G$ induces a group action on $V_\pm$. 
The properties  ($\G$1') -  ($\G$3')
follow directly from the corresponding properties 
($\G$1) -  ($\G$3) and Lemma \ref{Lem:orbitstabfinite}. 
We note that for Bernoulli percolation on $G = (V,E)$
the cluster decomposition of $V$ induces 
a random partition of $V_\pm := V_+ \cup V_-$.
Its distribution $\P$ is $\ph$-invariant for every $\ph \in \G$, 
because the distribution of Bernoulli percolation is invariant 
under all graph automorphisms. 
We have thus verified all assumptions of Theorem \ref{Thm:partitions}, 
and therefore we get the stochastic domination as claimed.
From this we immediately get the inequality for conditional expectations, 
since any $\{0,1,2,...\}$-valued random variable $\cX$ we have $\E(\cX) = \sum_{n \ge 0} \P(\cX \ge n)$. 
For the inequality of the unconditional expectations
we first note that in case of $\E_p(|\cC_\pm|) = \infty$
Lemma \ref{Lem:infinite} implies that $\E_p(|\cC_+|) = \infty = \E_p(|\cC_-|)$, since $(\G 1), (\G 2)$ imply the existence of some $\tau \in \G$ 
with $\tau(V_+) = V_-$ and $\tau(V_-) = V_+$. 
In case of $\E_p(|\cC_\pm|) < \infty$  we have 
$|\cC_\pm| < \infty$ $\P_p$-a.s., 
i.e. conditional and unconditional expectations are equal. 
\qed 

\bigskip

{\bf Proof of Proposition \ref{Prop:orbstab3}:}
In case of ($\G$3a) we are done by part (a) of Lemm \ref{Lem:orbstab}, 
and in case of  ($\G$3b) we adopt the argument for (b). 
Let $v,w \in V_\pm$ and set $k := d(v,w)$.  
Let $K \subset V$ be finite with $V_\pm \cap K \neq \emptyset$. 
Using transitivity, for every $x \in V_\pm \cap K$ we can choose 
$\ph_x \in \G$ such that $\ph_x(v) = x$ and set $x' := \ph_x(w)$. 
We note that $\ph_x(\G_v(w)) = \G_x(x')$, 
which implies $|\G_v(w)| = |\G_x(x')|$, 
and we note that $(x,y) \in \G(v,w) \Lra (x,y) \in \G(x,x') \Lra  y \in \G_x(x')$. This gives  
\begin{align*}
&\G(v,w) \cap K^2 
= \{(x,y) \in K^2: x \in V_\pm, y \in \G_x(x')\},   
\quad \text{ and thus }\\
&|V_\pm \cap K| |\G_v(w)| - |\G(v,w) \cap K^2|
= \sum_{x \in V_\pm \cap K} (|\G_x(x')| -|\G_x(x') \cap K| ).
\end{align*}
The last sum is $\ge 0$, and noting that for all $u \in \G_x(x')$ 
we have $d(x,u) = d(x,x') = d(v,w) = k$, we see that 
$\G_x(x') \subset K$ for $x \in V_\pm \cap (K \stm \partial_k K)$. Thus the above gives 
$$
0 \le |V_\pm \cap K| |\G_v(w)| - |\G(v,w) \cap K^2| \le \!\!\sum_{x \in V_\pm \cap \partial_k K} \! |\G_x(x')| 
 =  | V_\pm \cap \partial_k K| |\G_v(w)|. 
$$
By symmetry we obtain the same estimate with $v,w$ interchanged. 
Since  $|\G(v,w) \cap K| = |\G(w,v) \cap K^2|$
we see that 
\begin{align*}
&\Big||\G_v(w)| - |\G_w(v)|\Big| 
\le \frac{| V_\pm \cap \partial_k K| }{|V_\pm \cap K|} (|\G_v(w)| + |\G_w(v)|). 
\end{align*}
By ($\G$3b) the right hand side can be made arbitrarily small. \qed

\bigskip 

Finally we will infer the corollaries of Section \ref{Sec:results} from Theorem \ref{Thm:Graphs}. 
We simply need to  check conditions ($\G 1$)-($\G 3$) 
for a suitable choice  of $\G$. 
In fact, we always set 
$$
\G := \{\ph \in \Aut(G): \ph(V_+),\ph(V_-) \in \{V_+,V_-\}\},
$$ 
so that ($\G 1$) is trivially satisfied, 
we check  ($\G 2$) by describing suitable graph automorphisms, and we use 
the above proposition to check ($\G 3$). 

\bigskip 

{\bf Proof of Corollary \ref{Cor:bunkbed}:} 
For the given choice of $V_+ = V \times \{0\}$, $V_- = V \times \{1\}$ let $\G$ as above. 
$G \cp L_2$ is still connected and locally finite. By choice of $\G$ we have ($\G 1$).  
Making use of the the vertex-transitivity of $G$ and the reflection symmetry of $L_2$
we see that $G \cp L_2$ is vertex-transitive (and even vertex-swap-transitive if $G$ is so), 
and ($\G 2$) is satisfied. In order to check ($\G 3$) we use Proposition \ref{Prop:orbstab3}, 
and in case of  $*$-amenable $G$ we note that 
for finite $K \subset V$ we have 
$\partial_k (K \times \{0,1\}) = (\partial_k K) \times \{0,1\}$. \qed 

\bigskip

{\bf Proof of Corollary \ref{Cor:layer}:} This can be done analogously as the proof of Corollary \ref{Cor:bunkbed}. 
For any of the given choices of $V_+,V_-$ let $\G$ as above. 
In any case $V_+ = V \times A_+$, $V_- = V \times A_-$ for some $A_+,A_- \subset \Z$ and thus $V_\pm = V \times A_\pm$ 
for $A_\pm = A_+ \cup A_-$.
$G \cp L_\infty$ is still connected and locally finite. By choice of $\G$ we have ($\G 1$).  
We note that in all cases (a),(b) and (c) the symmetries of $A_\pm$ are acting swap-transitively on $A_\pm$. 
(For (a) consider the reflection in $\frac k 2$, for (b) and (c) consider the translations by $nl$ and the reflections in $nl+ \frac k 2$ for $l \in \Z$.)
This implies that ($\G 2$) is satisfied, and we obtain ($\G 3$) using Proposition \ref{Prop:orbstab3}. 
For this we note that swap-transitivity of $G$ implies that $\G$ acts swap-transitively on $V_\pm$, 
and that $*$-amenability of $G$ implies ($\G$3b) from Proposition \ref{Prop:orbstab3}. 
For the latter we consider finite  $K \subset V$ with $K \neq \emptyset$ and $I_N := \{-N,...,N-1\}$ and note that 
$\partial_k (K \times I_N) = (\partial_k K \times I_N) \cup  (K \times \partial_k I_N)$ and thus 
\begin{align*}
& \frac{|V_\pm \! \cap \partial_k (K \!\times\! I_N)|}{|V_\pm \! \cap (K \!\times\! I_N)|} \le 
\frac{|\partial_k K \!\times \! (A_\pm \!\cap I_N)|}{|K \!\times \! (A_\pm \!\cap I_N)|} + 
\frac{|K \!\times\! (A_\pm \!\cap \partial_k I_N)|}{|K \!\times \! (A_\pm \!\cap I_N)|}, 
\end{align*}
and on the right hand side the first gets arbitrarily small 
for an appropriate choice of $K$, 
and the second term tends to $0$ by choice of $A_\pm$. 
(In case (a) $A_\pm \!\cap \partial_k I_N = \emptyset$ for large $N$ and in cases (b),(c) 
$|A_\pm \!\cap \partial_k I_N| \le 2k$, whereas $|A_\pm \!\cap I_N| \to \infty$.) 
 \qed 

\bigskip 

{\bf Proof of Corollary \ref{Cor:Z2}:} This can be done analogously. 
For any of the given choices of $V_+,V_-$ let $\G$ as above. 
($\G 1$) holds and it is easy to check ($\G 2$): In case of (a) $\G$ contains all translations by $nv, n \in \Z$
and all reflections in the points $\frac n 2 v + \frac 1 2 w$, $n \in \Z$. 
In case of (b) $\G$ contains all translations by $nv + mu$, $n,m \in \Z$ and all reflections in the points 
$\frac n 2 v + \frac m 2 u + \frac 1 2 w$, $n,m \in \Z$. 
($\G 3$) follows via condition ($\G3b$) from Proposition \ref{Prop:orbstab3} choosing $K_N = \{-N,....,N\}^2$. 
In case of (a) $|V_\pm \cap K_N|$ grows linearly in $N$ and  $|V_\pm \cap \partial_k K_N|$ remains bounded. 
In case of (b) $|V_\pm \cap K_N|$ grows quadratic in $N$ and $|V_\pm \cap \partial_k K_N|$ grows linearly. 
 \qed 

\bigskip 

{\bf Proof of Corollary \ref{Cor:hypercube}:} 
We first note that $\P_p(o \lra v)$ indeed only depends on the 
graph distance $d(v,o)$, since for every permutation $\si$ of $\{1,...,d\}$ the map $\ph_\si$ define by $f_\si(v_1,...,v_d)= (v_{\si(1)},...,v_{\si(d)})$ is a graph automorphism of $L_2^d$ 
with $f_\si(o) = o$.
Let $k,l \ge 0$ such that $k+l \le d$.
W.l.o.g. we may assume $k \ge 1$, since otherwise the assertions are trivial. 
To obtain the first relation we consider the vertex sets 
$V_\pm = \{v \in V: \forall i > k+l: v_i = 0\}$, 
$V_+ = \{v \in V_\pm: \sum_{i=1}^k v_i \text{ is even}\}$ 
and  $V_-= \{v \in V_\pm: \sum_{i=1}^k v_i \text{ is odd}\}$. 
Let $\G := \{\ph \in \Aut(G): \ph(V_+),\ph(V_-) \in \{V_+,V_-\}\}$. 
Let $\tau_i: V \to V$ denote the reflection of the $i$-th component, 
i.e. $\tau_i(v) = w$ iff $w_j = v_j$ for all $j \neq i$ and $w_i = 1-v_i$. 
For all $i \le k+l$ we have $\tau_i \in \G$. 
This implies ($\G 2$). ($\G 3$) is satisfied by Proposition \ref{Prop:orbstab3}: $K = V$ is finite and $\partial_k K = \emptyset$. 
Thus Theorem \ref{Thm:Graphs} implies $\E(\cC_+) - \E(\cC_-) \ge 0$. 
Letting $A_{i,j} = \{v \in V_\pm: \sum_{m=1}^k v_m = i, \sum_{m = k+1}^{k+l} v_m = j\}$ we have 
\begin{align*}
&\E_p(\cC_+) - \E_p(\cC_-) 
= \sum_{v \in V_+} \P_p(o \lra v) - 
\sum_{v \in V_-} \P_p(o \lra v) \\
&= 
\sum_{i,j } \Big(\sum_{v \in V_+ \cap A_{i,j}} \P_p(o \lra v)
- \sum_{v \in V_- \cap A_{i,j}} \P_p(o \lra v)\Big)\\
&= \sum_{i,j } (-1)^i c_{i+j} |A_{i,j}| 
= \sum_{i,j } (-1)^i c_{i+j} \binom k i \binom l j. 
\end{align*}
This gives the first relation. For the second relation 
we choose $V_\pm := V_0 \cup V_1$, where 
$V_0 :=  \{v \in V: \forall i > k: v_i = 0\}$ and  
$V_1 :=  \{v \in V: \forall k < i \le k+l: v_i = 1, \forall i > k+l: v_i = 0\}$. We first set
$V_+ = \{v \in V_0: \sum_{i=1}^k v_i \text{ is even}\} 
\cup \{v \in V_1: \sum_{i=1}^k v_i \text{ is odd}\}$ 
and $V_- = \{v \in V_0: \sum_{i=1}^k v_i \text{ is odd}\} 
\cup \{v \in V_1: \sum_{i=1}^k v_i \text{ is even}\}$. 
Let $\G := \{\ph \in \Aut(G): \ph(V_+),\ph(V_-) \in \{V_+,V_-\}\}$. Let $\tau: V \to V$ denote the reflection of the components 
$k+1,...,k+l$, i.e. $\tau(v) = w$ iff $w_i = 1-v_i$ for all 
$i \in \{k+1,...,k+l\}$ and 
$w_i = v_i$ for all $i \notin \{k+1,...,k+l\}$. 
We have $\tau \in \G$ and $\tau_i \in \G$ for all $i \le k$. 
This implies ($\G 2$). 
($\G 3$) is satisfied by the same reasoning as above. 
Thus Theorem \ref{Thm:Graphs} implies $\E(\cC_+) - \E(\cC_-) \ge 0$. 
Letting $A_{i} = \{v \in V_\pm: \sum_{m=1}^k v_m = i\}$ we have 
\begin{align*}
&\E_p(\cC_+) - \E_p(\cC_-) 
= 
\sum_{i} \Big(\sum_{v \in V_+ \cap A_{i}} \P_p(o \lra v)
- \sum_{v \in V_- \cap A_{i}} \P_p(o \lra v)\Big)\\
&= \sum_{i} \Big((-1)^i c_{i} + (-1)^{i+1} c_{i+l}\Big)|A_{i}| 
= \sum_{i} (-1)^i c_{i} \binom k i  - \sum_{i} (-1)^i c_{i+l} \binom k i. 
\end{align*}
Similarly choosing $V_+ = \{v \in V_0: \sum_{i=1}^k v_i \text{ is even}\} 
\cup \{v \in V_1: \sum_{i=1}^k v_i \text{ is even}\}$ 
and $V_- = \{v \in V_0: \sum_{i=1}^k v_i \text{ is odd}\} 
\cup \{v \in V_1: \sum_{i=1}^k v_i \text{ is odd}\}$
we get 
\begin{align*}
0& \le \E_p(\cC_+) - \E_p(\cC_-) 
= \sum_{i} (-1)^i c_{i} \binom k i  + \sum_{i} (-1)^i c_{i+l} \binom k i. 
\end{align*}
In combination this gives the second relation.
Noting that $\sum_{i=0}^k (-1)^{k-i} \binom k i c_{m+i}$ 
$= \De^k[c](m)$ for all $k \le d$ and $m \le d-k$
we get the assertions for the discrete derivatives 
(in case of $l = 0$). 
\qed

\end{document}